# ON MINIMUM $(K_{1,r};k)$-VERTEX STABLE GRAPHS ON THE EXACT NUMBER OF VERTICES

ARTUR KUŹNAR

ABSTRACT. A graph $G$ is said to be $(H;k)$-vertex *stable* if $G$ contains a subgraph isomorphic to $H$ even after removing any $k$ of its vertices alongside with their incident edges. We will denote by $\mathrm{stab}(H;k)$ the minimum size among sizes of all $(H;k)$-vertex stable graphs. In this paper we consider a case where the structure $H$ is a star graph $K_{1,r}$ and the the number of vertices in $G$ is exact, i.e. equal to $1+r+k$. We will show that under the above assumptions $\mathrm{stab}(K_{1,r};k)$ equals either $\frac{1}{2}(k+1)(2r+k)$, $\frac{1}{2}\big((r+k)^2-1\big)$ or $\frac{1}{2}(r+k)^2$. Moreover, we will characterize all the extremal graphs.

## Contents

## 1. INTRODUCTION

The notion of $(H;k)$-vertex stability of a graph (also named $(H;k)$-vertex fault tolerance) stems directly from its applications to real-life problems. Consider the following challenge: given a net represented by a graph $H$, we want to immunize it against up to $k$ faults, while maintaining as low a cost of the whole structure as possible. We will work under two assumptions: the number of 'main' components, modeled as vertices of $H$, has to be as the minimum one (due to disparity of costs between them and the connections modeled as edges of $H$) and $H$ is a star graph that may be interpreted as an important junction in the net. Thus we only have to find the immune structure with the fewest number of connections. Observe, that both main components and connections may be faulty, the loss of a vertex is, however, much more detrimental, as it invalidates all connections the lost vertex had with the rest of the structure. Let us now formalize our approach.

In this paper we will only consider undirected graphs with no loops or multiple edges. Our notational system will closely mirror the one used by Diestel in [3]. Given a graph $G$, we denote its set of vertices by $V(G)$ and its set of edges by $E(G)$. Furthermore, $|G| := |V(G)|$ is the *order* of $G$ and $||G|| := |E(G)|$ is its *size*. The set of all neighbours of a given vertex $v$ will be denoted by $N(v)$. A vertex is considered *total* in $G$, if its set of neighbours is $V(G) \setminus \{v\}$. If $W \subseteq V$, by $G[W]$ we mean an induced subgraph of $G$ on the vertex set $W$. Unless it causes confusion, we will denote an edge with endpoints $u$ and $v$ simply as $uv$ instead of $\{u,v\}$. Similarly, if the vertices of a graph $G$ of order $n$ are labelled with natural numbers $1,\ldots,n$, we will denote the edges by $ij$ for some $i,j \in \{1,\ldots,n\}$. Since the set of initial $n$ positive integers will be used quite frequently, it is convenient to introduce a denotation $[n] := \{1,\ldots,n\}$. Given graphs $G_1$ and $G_2$, we define their *conjunction* $G_1 * G_2$ in the following way

$$V(G_1 * G_2) = V(G_1) \cup V(G_2),$$
$$E(G_1 * G_2) = E(G_1) \cup E(G_2) \cup E(G_1, G_2),$$

where $E(G_1, G_2)$ is the set of all edges between $V_1 = V(G_1)$ and $V_2 = V(G_2)$. When $v \in V$ and $W \subset V$ by $G - v$ and $G - W$ we mean respectively induced subgraphs $G[V \setminus \{v\}]$ and $G[V \setminus W]$ of $G$. A *star* graph is a complete bipartite graph $K_{1,r}$ for $r \in \mathbb{N}, r \geq 3$. By $K_n^{n-2}$ we mean an $(n-2)$-regular graph of order $n$ ($n$ is even).

A graph $G$ is said to be $(H,k)$-vertex *stable* if it contains a subgraph isomorphic to $H$ even after any $k$ of its vertices are removed with their incident edges. By $\text{stab}(G)$ we will denote the minimum size among the sizes of all $(H,k)$-vertex stable graphs. Observe that as long as $H$ has no isolated vertices, adding to or removing any number of vertices from an $(H;k)$-vertex stable graph produces a graph that is still $(H;k)$-vertex stable.

Since the introduction of the concept of vertex stability, many classes of graphs have been discussed. Cycles were discussed in [2] and [4], complete bipartite graphs in [5] and [6], with developments on complete graphs being the most advanced: [9], [10], [4], [8], [7]. We focus our attention on a class of star graphs $K_{1,r}$ for $r \geq 3$ and prove, that if the number of vertices is the smallest one necessary or *exact*, then $\text{stab}(K_{1,r};k)$ may assume one of the following values: $\frac{1}{2}(k+1)(2r+k)$, $\frac{1}{2}((r+k)^2 - 1)$ or $\frac{1}{2}(r+k)^2$, depending on the choice of $r$ and $k$. Moreover, all of the extremal graphs are characterized.

## 2. Construction of an $(H;k)$-vertex stable graph

In [1] Bruck, Cypher and Ho introduced a simple way of constructing an $(H;k)$-vertex stable graph on the exact number of vertices. Let $H$ be a given graph of order $n$ and let $k$ be a non-negative integer. We proceed as follows:

(1) Label all the vertices of $H$ with different elements of $[n]$,
(2) Add $k$ vertices to $H$, labelled with consecutive integers $n+1,\ldots,n+k$,
(3) if $ij \in E(H)$, connect each vertex in $\{i,\ldots,i+k\}$ with each one in $\{j,\ldots,j+k\}$ (unless they already are connected).

Recall the reasoning behind the $(H;k)$-vertex stability of a graph obtained this way. Let $F$ be a $k$-element set of faults. Consider a mapping

$$\psi : V(H) \to V(G - F)$$

assigning each vertex $i \in V(H)$ with the smallest unassigned label among labels of vertices in $V(G - F)$. Clearly $|V(H)| = |V(G - F)|$. Furthermore, if $ij \in E(H)$, then $\psi(i)\psi(j) \in E(G - F)$, since

$$\psi(i) \in \{i,\ldots,i+k\} \text{ and } \psi(j) \in \{j,\ldots,j+k\}$$

and the construction guarantees the existence of all the edges between the above sets. Thus $G - F$ contains a subgraph isomorphic to $H$. Note, that although each vertex is assigned a different value from $[n+k]$, the rule according to which the vertices of $H$ are labelled is not specified. As the procedure depends on the labelling, the resulting $(H;k)$-vertex stable graph is not uniquely determined if only $H$ is given, as may be seen in the following example.

**Example 1.** Let $H$ be a graph with $V(H) = \{u,v,w,t\}$ and $E(H) = \{uv, uw, ut, wt\}$. Consider the following labellings of vertices of $H$:

$$\eta_i : V(H) \to \{1,2,3,4\}, \text{ for } i = 1,2.$$

We define them as follows:
$$\eta_1(u) = 3, \ \eta_1(v) = 4, \ \eta_1(w) = 1, \ \eta_1(t) = 2$$
$$\eta_2(u) = 1, \ \eta_2(v) = 2, \ \eta_2(w) = 3, \ \eta_2(t) = 4.$$

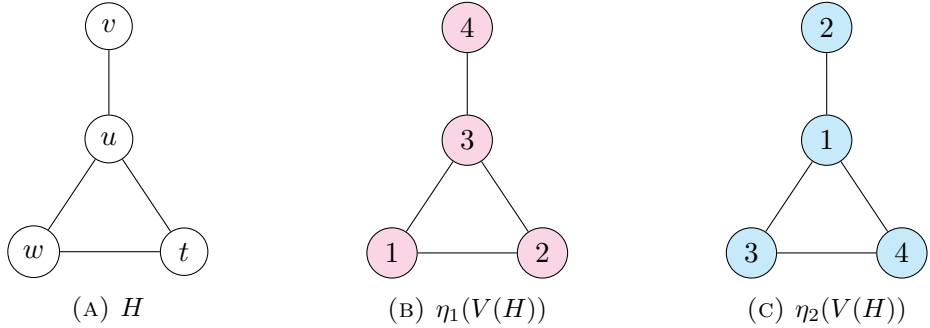

FIGURE 1. Labellings of vertices of $H$

(A) $H$  (B) $\eta_1(V(H))$  (C) $\eta_2(V(H))$

By the above construction, two non-isomorphic $(H;2)$-vertex stable graphs $G_1$ and $G_2$ are obtained. In both cases $V(G_1) = V(G_2) = [6]$, however

$E(G_1) = \{\{1,2\},\{1,3\},\{1,4\},\{1,5\},\{2,3\},\{2,4\},\{2,5\},\{3,4\},\{3,5\},\{3,6\},\{4,5\},\{4,6\},\{5,6\}\}$,

$E(G_2) = \{\{i,j\} : i,j \in \{1,2,\ldots,6\} \land i \neq j\}$.

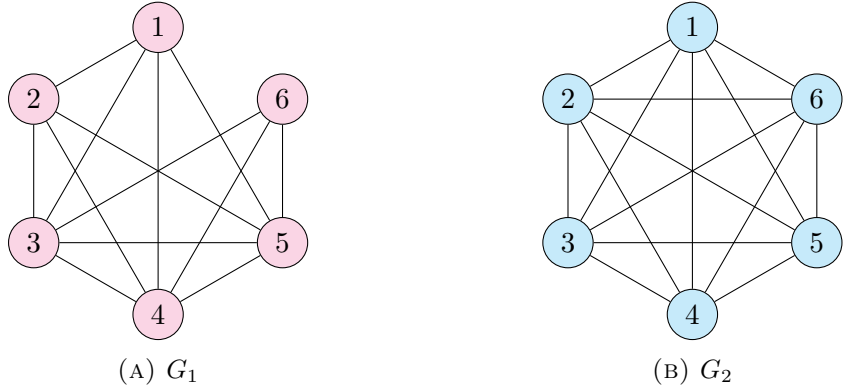

FIGURE 3. $(H;2)$-vertex stable graphs obtained respectively for $\eta_1$ and $\eta_2$

(A) $G_1$  (B) $G_2$

Observe, that the procedure may yield $(H;k)$-vertex stable graphs that are not only non-isomorphic but differ in size as well. Given a graph $H$ and labelling $\eta$ of vertices of $H$, denote by $G(H;k;\eta)$ the $(H;k)$-vertex stable graph obtained via Bruck, Cypher and Ho's method. Our goal will be now to show, that in specific cases of our interest, the construction produces a unique graph, with no regard for the labelling.

**Lemma 1.** *Let $K_{1,r}$ be a star graph. Then for each non-negative integer $k$ a $(K_{1,r};k)$-vertex stable graph obtained via from Bruck, Cypher and Ho's construction is determined uniquely. Moreover, it is isomorphic with $K_{k+1} * \overline{K}_r$.*

*Proof.* We will show first, that the choice of label for the central vertex of the star graph does not affect the size of $G(K_{1,r};k;\eta)$. Let $\eta$ be a labelling that assigns the central vertex some

$t \in [r+1]$. Clearly, the leaves are assigned values from $[r+1] \setminus \{t\}$, thus
$$E(K_{1,r}) = \bigcup_{l \in [r+1] \setminus \{t\}} \{\{l,t\}\},$$
$$V(G(K_{1,r};k;\eta)) = \{1, 2, \ldots, r+1, r+2, \ldots, \ldots, r+k+1\}.$$

Our next step is to determine the size of $G(K_{1,r};k;\eta)$. By definition of $G(K_{1,r};k;\eta)$ we obtain:

$$\begin{aligned} E(G(K_{1,r};k;\eta)) &= \bigcup_{l \in [r+1] \setminus \{t\}} \{ij : i \in \{t, \ldots, t+k\} \wedge j \in \{l, \ldots, l+k\} \wedge i \neq j\} \\ &= \Big\{ij : i \in \{t, \ldots, t+k\} \wedge j \in \bigcup_{l \in [r+1] \setminus \{t\}} \{l, \ldots, l+k\} \wedge i \neq j\Big\} \\ &= \Big\{ij : i \in \{t, \ldots, t+k\} \wedge j \in [r+k+1] \wedge i \neq j\Big\}. \end{aligned}$$

Observe, that $\bigcup_{l \in [r+1] \setminus \{t\}} \{l, \ldots, l+k\} = [r+k+1]$ if and only if $t \notin \{1, r+1\}$. If $t \in \{1, r+1\}$, we have
$$\bigcup_{l \in [r+1] \setminus \{t\}} \{l, l+1, \ldots, l+k\} = [r+k+1] \setminus \{t\}.$$

Note, however, that even when $t \in \{1, r+1\}$ for every vertex in $\{t, \ldots, t+k\}$ there exists an edge that connects it with each other vertex with the exception of $t$. On the other hand, $i$ assumes all the values from $\{t, \ldots, t+k\}$, which guarantees the existence of all edges between $t$ and $[r+k+1] \setminus \{t\}$. Therefore
$$\Big\{ij : i \in \{t, \ldots, t+k\} \wedge j \in [r+k+1] \wedge i \neq j\Big\} = \Big\{ij : i \in \{t, \ldots, t+k\} \wedge j \in [r+k+1] \setminus \{t\} \wedge i \neq j\Big\}.$$

Consider next $i', j' \in [r+k+1] \setminus \{t, t+1, \ldots, t+k\} = \{1, \ldots, t-1\} \cup \{t+k+1, \ldots, r+k+1\}$. By the construction of $(G(K_{1,r};k;\eta)$, there have to exist such $i, j \in V(K_{1,r})$ that $i' = \psi(i)$ and $j' = \psi(j)$. Hence $i' \in \{i, \ldots, i+k\}$ and $j' \in \{j, \ldots, j+k\}$. Were $i'j'$ an edge of $G(K_{1,r};k;\eta)$, there would have be either $i = t$ or $j = t$. It is, however, impossible for otherwise one of $i', j'$ would have to be in $\{t, \ldots, t+k\}$. Therefore $\{1, \ldots, t-1\} \cup \{t+k+1, \ldots, r+k+1\}$ is an independent set. We have seen so far, that

(1) For $i \in \{t, t+1, \ldots, t+k\}$
$$N(i) = [r+k+1] \setminus \{i\}, \text{ hence } |N(i)| = r+k.$$

(2) For $i \in \{1, \ldots, t-1\} \cup \{t+k+1, \ldots, r+k+1\}$ we have only connections with the remaining $k+1$ vertices. Thus $|N(i)| = k+1$.

By Euler's handshake lemma we immediately get the size of $G(K_{1,r};k;\eta)$:
$$|E(G(K_{1,r};k;\eta))| = \frac{1}{2} \sum_{v \in V(G(K_{1,r};k;\eta))} \deg(v) = \frac{1}{2}((k+1)(r+k) + r(k+1)) = \frac{1}{2}(k+1)(2r+k).$$

We have thus proven, that the size of $G(K_{1,r};k;\eta)$ is independent of the chosen label for the central vertex of $K_{1,r}$. Owing to the symmetry of a star graph, it is also independent of the labels of leaves – and as a consequence does not depend on the labelling $\eta$ of $K_{1,r}$ at all.

It remains to be seen, that $G(K_{1,r};k;\eta)$ is uniquely determined, up to isomorphism. Let $G(K_{1,r};k;\eta_1)$ and $G(K_{1,r};k;\eta_2)$ be two $(K_{1,r};k)$-vertex stable graphs obtained via Bruck, Cypher and Ho's procedure using two distinct labellings $\eta_1$ and $\eta_2$. Take any bijective mapping $\zeta : V(G_1) \to V(G_2)$, satisfying the following conditions:
$$\zeta(\{t_1, \ldots, t_1+k\}) = \{t_2, \ldots, t_2+k\}$$
and
$$\zeta(\{1, \ldots, t_1-1\} \cup \{t_1+k+1, \ldots, r+k+1\}) = \{1, \ldots, t_2-1\} \cup \{t_2+k+1, \ldots, r+k+1\}.$$

Let $i,j \in V(G_1)$ and let $ij \in E(G_1)$. Then either $i \in \{t_1,\ldots,t_1+k\}$ or $j \in \{t_1,\ldots,t_1+k\}$. By the definition of $\zeta$, we obtain that either $\zeta(i) \in \{t_2,\ldots,t_2+k\}$ or $\zeta(j) \in \{t_2,\ldots,t_2+k\}$, hence $\zeta(i)\zeta(j) \in E(G_2)$.

Conversely, let now $i',j' \in V(G_2)$ and let $i'j' \in E(G_2)$. Then either $i' \in \{t_2,\ldots,t_2+k\}$ or $j' \in \{t_2,\ldots,t_2+k\}$, which implies that either $\zeta^{-1}(i') \in \{t_1,\ldots,t_1+k\}$ or $\zeta^{-1}(j') \in \{t_1,\ldots,t_1+k\}$, thus $\zeta^{-1}(i')\zeta^{-1}(j') \in E(G_1)$.

It is clear now, that the construction yields for $H = K_{1,r}$ only one graph (up to isomorphism), which we will denote by $G(r,k)$ for the sake of simplicity. $\square$

FIGURE 5. $K_8 * \overline{K}_8$ – a $(K_{1,8};7)$-vertex stable graph obtained via Bruck, Cypher and Ho's construction.

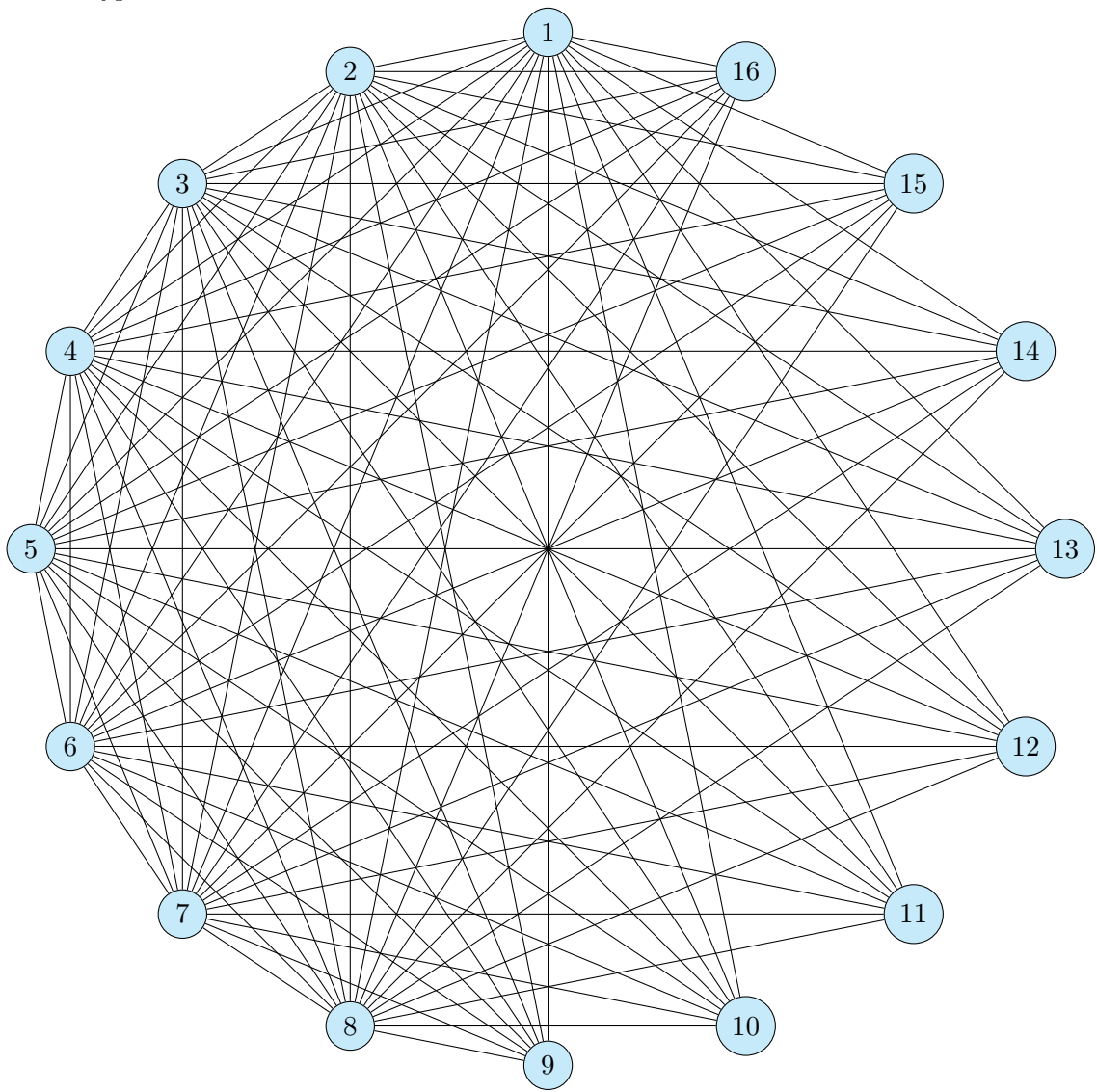

## 3. Minimum $(K_{1,r};k)$-vertex stable graphs

In the previous section we have seen that there is a non-complicated way of constructing $(K_{1,r};k)$-vertex stable graphs, such, that the resulting graph is uniquely determined. A question arises, whether or not the graph $G(r,k)$ is the one of the smallest size among all $(K_{1,r};k)$-vertex stable graphs. Here is the answer:

**Theorem 1.** *Let $G$ be a $(K_{1,r};k)$-vertex stable graph of order $r+k+1$ and let $k_0 = \min\{k' : k' > (r-1)^2 - 2 \wedge 2 \nmid k'\}$.*

*(1) Assume at least one of the following conditions holds:*
  *a) $r$ is odd,*
  *b) $r$ is even and $k < (r-1)^2 - 2$,*
  *then*
  $$|E(G)| \geq \frac{1}{2}(k+1)(2r+k)$$
  *and the sole extremal graph is $G(r,k)$.*

*(2) If $r$ is even and $k = (r-1)^2 - 2$ (respectively $k = (r-1)^2 - 1$), then*
  $$|E(G)| \geq \frac{1}{2}(k+1)(2r+k)$$
  *and there are exactly two extremal graphs: $G(r,k)$ and $K_{r+k+1}^{r+k-1}$ (respectively $G(r,k)$ and $K_{r+k}^{r+k-2} * K_1$).*

*(3) If $r$ is even, $k$ is odd and satisfies $k \geq k_0$, then*
  $$|E(G)| \geq \frac{1}{2}((r+k)^2 - 1)$$
  *and the sole extremal graph is $K_{r+k+1}^{r+k-1}$.*

*(4) If both $r$ and $k$ are even and $k$ satisfies $k \geq k_0$, then*
  $$|E(G)| \geq \frac{1}{2}(r+k)^2$$
  *and the sole extremal graph is $K_{r+k}^{r+k-2} * K_1$.*

To prove this theorem, we apply the following useful lemma.

**Lemma 2.** *The only $(K_{1,r};k)$-vertex stable graph $G$ of order $r+k+1$ and such, that $\Delta(G) < r+k$, $r$ is even and $k$ is odd, is $K_{r+k+1}^{r+k-1}$. Moreover, if $r$ is odd or $k$ is even, no graph of order $r+k+1$ with $\Delta < r+k$ is $(K_{1,r};k)$-vertex stable.*

*Proof of the lemma.* Let $G$ be a $(K_{1,r};k)$-vertex stable graph of order $r+k+1$ with $\Delta(G) < r+k$. It is clear, that for every vertex $v \in V(G)$ we can find such $u \in V(G)$, that $uv \notin E(G)$. It is convenient to view the process of deleting $k$ vertices from $G$ as introducing a partition of $V(G)$ into two classes, the elements of which may be moved from one to the other. Let $P = (S, F)$ be such a partition. The set $S$ is the set of vertices remaining in $G$, on which we hope to find a star graph. The set $F$ is on the other hand a set of faults, i.e. the set of deleted vertices. Let $\mathcal{P}$ denote the family of all such partitions $P = (S, F)$, that $|S| = r+1$ and $|F| = k$. If in the partition $P \in \mathcal{P}$ for every vertex $v \in S$ there existed such $u \in S$, that $uv \notin E(G)$, then $\Delta(G[S]) \leq r-1$ and $G$ would certainly not be $(K_{1,r};k)$-vertex stable. Let now $P = (S, F)$ be a partition, for which the number $s$ of vertices of degree $r$ in $S$ is the smallest (though non-zero) among all partitions $P' \in \mathcal{P}$. Suppose first, that $s \geq 2$; then there exists $\{w_1, w_2, \ldots, w_s\} = W \subseteq S$ satisfying
$$\deg_{G[S]}(w_1) = \deg_{G[S]}(w_2) = \ldots = \deg_{G[S]}(w_s) = r.$$

As $\Delta(G) < r+k$, there exists as well such a set $\{u_1, u_2, \ldots, u_t\} = U \subseteq F$, that for every vertex $w \in W$ we can find $u \in U$ with $uw \notin E(G)$. If there are $w_1, w_2 \in W$, $u_1, u_2 \in U$ :

$u_i w_i \notin E(G)$ ($i = 1, 2$), consider a partition $P' = (S', F')$:
$$S' = S \setminus \{w_2\} \cup \{u_1\},$$
$$F' = F \setminus \{u_1\} \cup \{w_2\}.$$

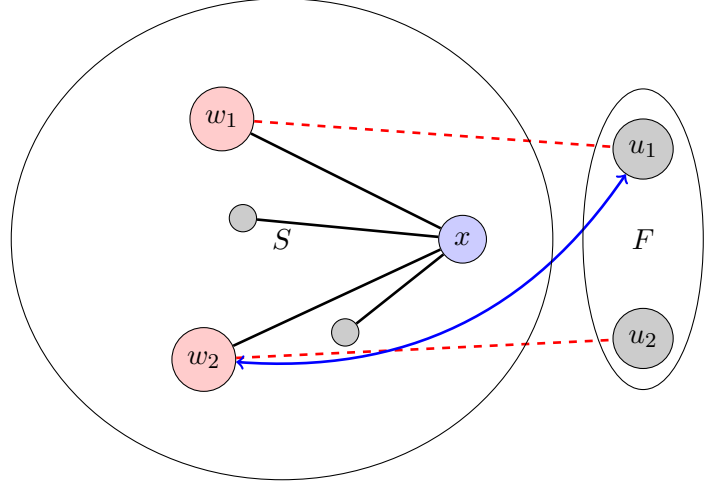

FIGURE 6. In $S$ there are at least two vertices $w_1, w_2$ of degree $r$ in $G[S]$.

In the above auxiliary illustration the red dotted lines are used to represent edges of $\overline{G}$, while the blue arrows mark the exchange. We will make use of this convention while considering the next cases as well.

Observe, that $\deg_{G[S]}(v) < r$ implies $\deg_{G[S']}(v) \leq \deg_{G[S]}$. Indeed, removing a vertex $w_2$ of degree $r$ in $G[S]$ decreases the degrees of all other vertices in $S$ by 1. Adding $u_1$ increases the degrees of those vertices by at most 1, so for all vertices $v \in S \setminus \{w_2\}$ we have
$$\deg_{G[S]}(v) - 1 \leq \deg_{G[S']}(v) \leq \deg_{G[S]}(v);$$
therefore none of the vertices in $S \setminus W$ has degree $r$ in $G[S']$. Thus
$$|\{v \in S' : \deg_{G[S']}(v) = r\}| = s' \leq s - 2.$$

When vertices $u_1$ and $u_2$ are not distinct, we may construct $S'$ by deleting any of the vertices $w_1, w_2$ and adding $u_1 = u_2$. Similarly, no vertex of degree less than $r$ in $G[S]$ will have degree $r$ in $G[S']$, however
$$s' = |W'| = |\{w \in S' : \deg_{G[S']}(w) = r\}|$$
will decrease by at least 2. In both cases we come to a contradiction with the choice of $P$. Furthermore, we are able to reduce the number of vertices of degree $r$ to zero or one in a finite number of steps. As $s = 0$ implies that $G$ is not ($K_{1,r}; k$)-vertex stable, we may consider only the case when $s = 1$. We may single out the vertex $x \in S : \deg_{G[S]}(x) = r$. We will complete the proof in five steps.

**Step 1:** If there are such vertices $y, y_1, y_2$ in $S \setminus \{x\}$, that $yy_1 \notin E(G[S])$ and $yy_2 \notin E(G[S])$, then having in mind that there is in $F$ such a vertex $x'$ that, $xx' \notin E(G)$, we may construct partitions $P^{(i)} = (S^{(i)}, F^{(i)})$ in the following way:
$$S^{(i)} = S \setminus \{y_i\} \cup \{x'\},$$
$$F^{(i)} = F \setminus \{x'\} \cup \{y_i\},$$
where $i \in \{1, 2\}$.

FIGURE 7. There are two distinct vertices $y_1, y_2$ in $S$ not connected with $y$.

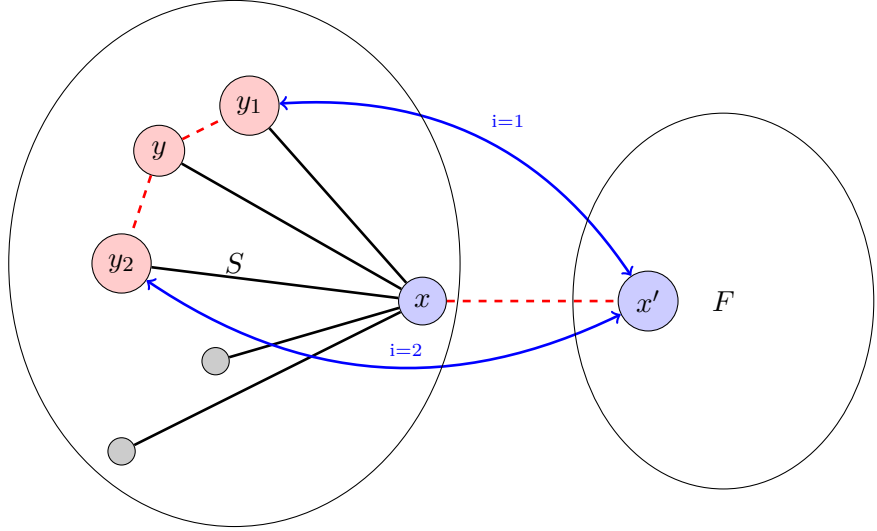

Observe, that again none of the vertices in $S \setminus \{x\}$ has degree $r$, though we cannot use the same reasoning to show it, owing to the degrees of $y_1, y_2$ and $y$ in $G[S]$. If for all the vertices in $S \setminus \{x, y, y_1, y_2\}$, all of their corresponding non-neighbouring (in $G$) vertices are in this set as well, the exchange is safe. Problems may arise if there are such vertices $y'_1, y'_2 \in S \setminus \{x\}$, for which respectively $y_1$ and $y_2$ are the only non-neighbouring vertices in $G[S \setminus \{x\}]$ and $x'$ is connected with both $y'_1$ and $y'_2$. We may, however, construct the partitions $P^i$ in the following way:

$$S' = S \setminus \{y'_i\} \cup \{x'\},$$
$$F' = F \setminus \{x'\} \cup \{y'_i\},$$

where $i \in \{1, 2\}$.

FIGURE 8. There are two distinct vertices $y_1, y_2$ in $S$ not connected with $y$ – second case.

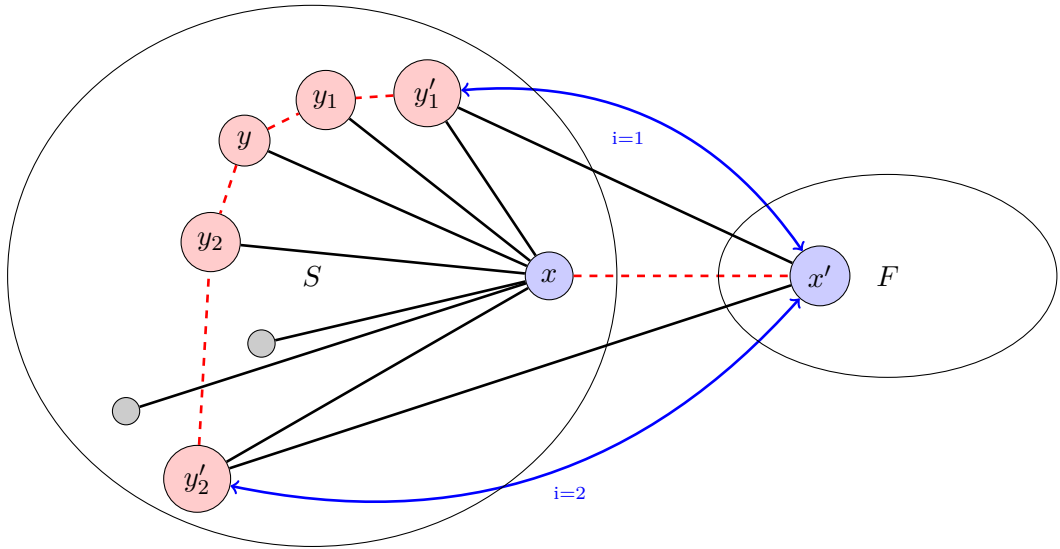

For both of $P' \in \{P^1, P^2\}$ we get $s' = 0$, hence the contradiction with the choice of $P$. Therefore for every vertex $v \in S \setminus \{x\}$ there exists exactly one such $u \in S \setminus \{x\}$, that $uv \notin E(G[S])$, which implies the even number of elements in $S \setminus \{x\}$.
**Step 2:** If there are vertices $x'$ and $x''$ in $F$ such, that $xx' \notin E(G)$ and $xx'' \notin E(G)$, we may create a partition $P' = (S', F')$:

$$S' = S \setminus \{u, v\} \cup \{x', x''\},$$
$$F' = F \setminus \{x', x''\} \cup \{u, v\},$$

where $u, v \in S$ and $uv \notin E(G[S])$. Then in $S'$ there are no vertices of degree $r$ in $G[S']$, which leads us to contradiction with the choice of $P$. Then there exists exactly one such $x' \in V(G)$, that $xx' \notin E(G)$, hence $\deg_G(x) = r + k - 1$.

FIGURE 9. There are two vertices in $F$, none of which is connected with $x$.

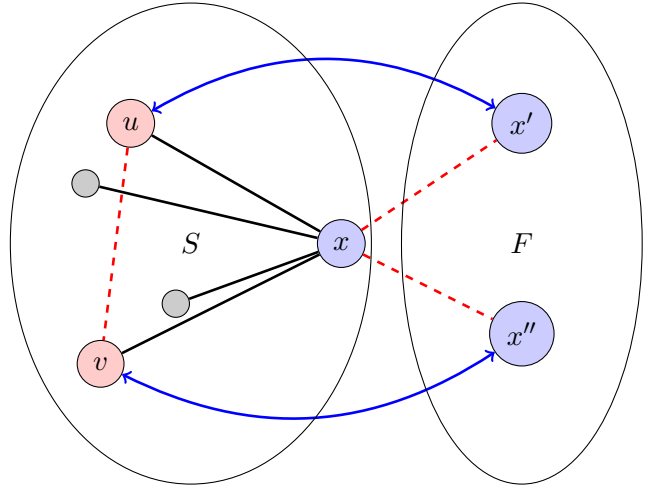

**Step 3:** Fix now the set $S \setminus \{x\}$. For all $y \in F \cup \{x\}$ create then a corresponding partition $P^{(y)} = \{S^{(y)}, F^{(y)}\}$ in the following way:

$$S^{(y)} = S \setminus \{x\} \cup \{y\},$$
$$F^{(y)} = F \setminus \{y\} \cup \{x\}.$$

Since $x$ was the only vertex of degree $r$ in $G[S]$, then $\deg_{G[S^{(y)}]}(y) = r$, for otherwise we would come to a contradiction with the choice of $P$. As $P^{(y)}$ is equivalent to $P$ with respect to the number of vertices of degree $r$ remaining respectively in sets $S^{(y)}$ and $S$, we may repeat the reasoning from previous steps for $P$. Therefore, for every $y \in F \cup \{x\}$ the degree $\deg_G(y) = r + k - 1$, which implies, that for every $y \in F \cup \{x\}$ there exists exactly one such $y' \in F \cup \{x\}$, that $yy' \notin E(G)$ and forces the even cardinality of $F \cup \{x\}$.
**Step 4:** Let once again $x' \in F$ be a vertex satisfying $xx' \notin E(G)$. By the above observations, there are exactly $\frac{|S \setminus \{x\}|}{2}$ pairs of vertices $\{u, u'\}$ in $S \setminus \{x\}$ such, that $uu' \notin E(G[S])$. Similarly, there are exactly $\frac{|F \setminus \{x'\}|}{2}$ pairs of vertices in $\{v, v'\}$ in $F \setminus \{x'\}$ such, that $vv' \notin E(G)$. Moreover, these numbers are integers and $\deg_G(v) = \deg_G(v') = r + k - 1$. Create a partition $P^{(u,v)} = (S^{(u,v)}, F^{(u,v)})$ as follows:

$$S^{(u,v)} = S \setminus \{u, u'\} \cup \{v, v'\},$$
$$F^{(u,v)} = S \setminus \{v, v'\} \cup \{u, u'\}.$$

FIGURE 10. We simultanously exchange pairs of vertices between $S$ and $F$.

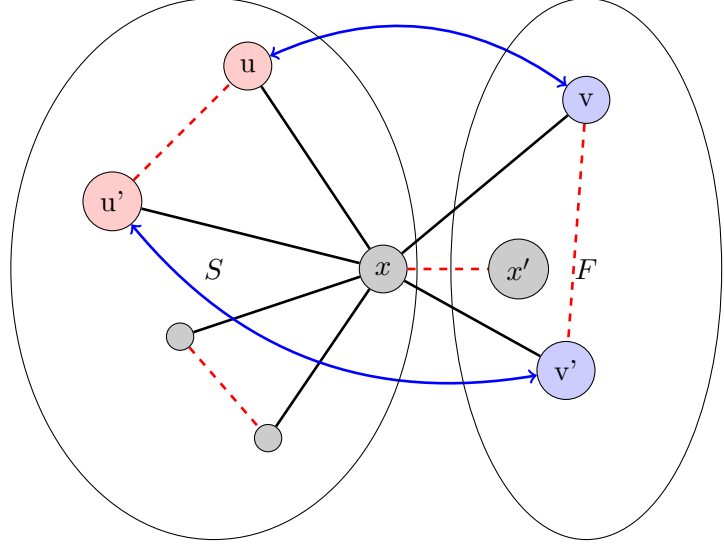

For every vertex $w \in S^{(u,v)} \setminus \{x\}$ there still exists a vertex $w'$ in $S^{(u,v)}$ such, that $ww' \notin E(G)$. Thus for the partition $P^{(u,v)}$, its corresponding number $s^{(u,v)} \leq 1$. In fact, $s^{(u,v)} = 1$, for otherwise we would come to a contradiction with the choice of $P$. We may next repeat the reasoning from step 3 for vertices $u, u' \in F^{(u,v)} \cup \{x\}$, creating partitions $P^{(u)}$ and $P^{(u')}$; note, that $\deg_G(u) = \deg_G(u') = r+k-1$. Since the choice of $\{u, u'\}$ is arbitrary, we conclude that all vertices $v \in V(G)$ have degree $\deg_G(v) = r+k-1$.

**Step 5:** It remains to be shown, that the graph $K_{r+k+1}^{r+k-1}$ is actually $(K_{1,r};k)$-vertex stable for even $r$ and odd $k$. Let $F$ be a $k$-element subset of $V(K_{r+k+1}^{r+k-1})$. Since the cardinality of $V(K_{r+k+1}^{r+k-1}) \setminus F$ is odd, there is at least one unpaired vertex in $K_{r+k+1}^{r+k-1} - F$ (its sole non-neighbouring vertex was deleted); it is clearly of degree $r$. Moreover, no other graph $G$ for which $\Delta(G) < r+k$ is $(K_{1,r};k)$-vertex stable. $\square$

*Proof of the theorem.* Having laid the groundwork, we may now proceed with the proof of the main theorem. We will apply inductive reasoning with cases building upon their predecessors.
**Case (1)** Let $G$ be a $(K_{1,r};k)$-vertex stable graph of order $r+k+1$, satisfying (1). When $k=0$, $G$ has exactly $r+1$ vertices, the degree of one of them being equal to $r$. Thus

$$|E(G)| \geq r = \frac{1}{2}(0+1)(2r+0),$$

the only $G$ for which the equality holds being simply $K_{1,r}$. Obviously, no other graph of size $r$ on $r+1$ vertices is $(K_{1,r}, 0)$-vertex stable. Suppose now, that the inequality holds true for all natural $k' < k$ and consider the following cases.

Let first $\Delta(G) \leq r+k-1$ for some $k$. By Lemma 2, $G$ is then either $(r+k-1)$-regular, or is not $(K_{1,r};k)$-vertex stable at all. Since the graphs that are not vertex stable are of no interest to us, we may assume $G$ is regular, $r$ is even and $k$ is odd. Suppose, that $||G|| < ||G(r,k)||$. Then

$$|E(G)| < |E(G(r,k))| \iff \frac{1}{2}(r+k+1)(r+k-1) \leq \frac{1}{2}(k+1)(2r+k)$$
$$\iff k^2 + 2kr + r^2 - 1 \leq k^2 + 2kr + 2r + k$$
$$\iff r^2 - 2r + 1 < k + 2$$
$$\iff (r-1)^2 - 2 \leq k.$$

Suppose next, that there is a vertex $v$ in $G$ with $\deg_G(v) = r + k$. Clearly, $G' = G - v$ obtained from $G$ by removing $v$ with all its incident edges is still $(K_{1,r}, k-1)$-vertex stable. Furthermore
$$|V(G')| = |V(G - v)| = |V(G)| - 1 = (k - 1) + r + 1.$$
By inductive hypothesis
$$|E(G')| \geq \frac{1}{2}(k - 1 + 1)(2r + k - 1) = \frac{1}{2}k \cdot (2r + k - 1),$$
from whence the lower bound on the size of $G$ may be immediately derived:
$$|E(G)| = |E(G')| + r + k \geq \frac{1}{2}(2kr + k^2 + 2k + 2r - k) = kr + \frac{1}{2}k^2 + \frac{1}{2}k + r$$
$$= \frac{1}{2}(2kr + k^2 + k + 2r) = \frac{1}{2}(k + 1)(2r + k) = ||G(r, k)||.$$
It remains to show, that the extremal graph is unique up to isomorphism. Observe, that there is a sequence of vertices $v_1, \ldots, v_k \in V(G)$ such, that for all $i \in [k]$
$$\deg_{G^{(i-1)}}(v_i) = r + k + 1 - i,$$
where $G^{(0)} = G$ and for all $i \in [k]$ graphs $G^{(i)}$ are defined inductively by $G^{(i)} = G^{(i-1)} - v_i$. Thus we get the following lower bound for the size of $G$:
$$|E(G)| = \sum_{i=1}^{k} \deg_{G^{(i-1)}}(v_i) + ||G^{(k)}|| = \sum_{l=1}^{k}(r + l) + ||G^{(k)}|| = kr + \frac{k(k+1)}{2} + ||G^{(k)}||$$
$$\geq kr + r + \frac{k(k+1)}{2} = \frac{1}{2}(k + 1)(2r + k),$$
with left hand side being equal to the right one if and only if $||G(k)|| = ||K_{1,r}||$ – which in turn takes place exactly when $G^{(k)} \cong K_{1,r}$. Denote the vertices of $K_{1,r}$ by $v_{k+1}, \ldots, v_{r+k+1}$, without the loss of generality setting
$$\deg_{K_{1,r}}(v_{k+1}) = r \ \wedge \ \deg_{K_{1,r}}(v_{k+2}) = \deg_{K_{1,r}}(v_{k+3}) = \ldots = \deg_{K_{1,r}}(v_{r+k+1}) = 1. \qquad (1)$$
Reconstructing inductively $G$ from $G^{(k)}$, we will now prove two claims:

**Claim 1.** For all $i \in [k] \cup \{0\}$, $\Delta(G^{(i)}[\{v_{k+2}, \ldots, v_{r+k+1}\}]) = 0$.

**Claim 2.** For all $i \in [k] \cup \{0\}$ and for all $j \in [k+1]$, $\deg_{G^{(i)}}(v_j) = r + k - i$.

When $i = k$, $G^{(i)}$ is a star graph and both Claim 1 and Claim 2 are trivially satisfied due to (1). Suppose now that they also hold true for all $i' > i$. Then $|G^{(i)}| = r + k + 1 - i$ and the graph $G^{(i)}$ is obtained through the addition of a vertex $v_{i+1}$ of degree $r + k - i$. Since the set $\{v_{k+2}, \ldots, v_{r+k+1}\}$ is independent in $G^{(i+1)}$ by inductive hypothesis and no edges are added between its elements, it remains such in $G^{(i)}$. Furthermore, also by inductive hypothesis, the degrees of all the vertices $v_{k+1}, v_k, \ldots, v_{i+2}$ are equal to $r + k - i - 1$ in $G(i + 1)$, thus adding a total vertex to $G(i+1)$ increases their degree by one and we get
$$\deg_{G^{(i)}}(v_{k+1}) = \ldots = \deg_{G^{(i)}}(v_{i+2}) = \deg_{G^{(i)}}(v_{i+1}) = r + k - i,$$
which concludes the proof of both claims. In particular, having continued the reconstructing process, we may formulate a new claim:

**Claim 3.** For all $i \in [k + 1]$, the degree $\deg_G(v_i) = r + k$, while for all $i \in [r + k + 1] \setminus [k + 1]$, the degree $\deg_G(v_i) = |\{v_1, \ldots, v_{k+1}\}| = k + 1$.

We have already proved in Lemma 1 that $G(r, k)$ is uniquely determined up to isomorphism; making use of a labelling with $t = 1$, we obtain that $\deg_{G(r,k)}(i) = r + k$ for all $i \in [k + 1]$ and $\deg_{G(r,k)}(i) = k + 1$ for $i \in [r + k + 1] \setminus [k + 1]$. It suffices to consider such a mapping $\zeta : V(G) \to V(G(r, k))$, that $\zeta(v_i) = i$ for all $i \in [r + k + 1]$. It is indeed an isomorphism, since $v_i v_j \in E(G)$ for some $i, j \in [r+k+1]$ if and only if $\zeta(v_i)\zeta(v_j) = ij \in E(G(r,k))$, for otherwise we would have
$$(i \in [k+1] \lor j \in [k+1]) \land i, j \in [r+k+1] \setminus [k+1],$$
a clear contradiction. Hence $G \cong G(r,k)$.

**Case (2)** We will now consider borderline cases, where the extremal graph is not uniquely determined. Note, that even value of $r$ implies that $(r-1)^2 - 2$ is odd. Let $r$ be even and let $k_1 := (r-1)^2 - 2$. Consider two $(K_{1,r}, k_1)$-vertex stable graphs $G_1$ and $G_2$. Suppose there exists a vertex $v$ in $G_1$ of degree $r + k_1$. Deleting it with all incident edges yields a graph $G'_1 = G_1 - v$, that is still $(K_{1,r}, k_1 - 1)$-vertex stable. However, $k_1 - 1 < (r-1)^2 - 2$ and $G'_1$ belongs to the class already considered in (1), therefore

$$|E(G_1)| = \deg_{G_1}(v) + |E(G'_1)| \geq (r + k_1) + \frac{1}{2}(k_1 - 1 + 1)(2r + k_1 - 1)$$
$$= (r + k_1) + \frac{1}{2}k_1(2r + k_1 - 1)$$
$$= \frac{1}{2}(k_1 + 1)(2r + k_1)$$

with both sides being equal if and only if $G'_1$ is extremal – thus isomorphic with $G(r, k_1 - 1)$. Suppose next $\Delta(G_2) < r + k_1$; by Lemma 2 the graph $G_2$ is then $(r + k_1 - 1)$-regular, hence

$$|E(G_1)| = \frac{1}{2}(r + k_1 + 1)(r + k_1 - 1) = \frac{1}{2}((r + k_1)^2 - 1).$$

It is clear, that $||G_1|| = ||G_2||$. As we have exhausted all the possibilities in this case, we conclude that $G_1$ and $G_2$ are the only extremal graphs.

Consider now a $(K_{1,r}; k_1 + 1)$-vertex stable graph $G$. Since $k_1$ is odd, $k_1 + 1$ is even; then by Lemma 2 there exists a vertex $v$ in $G$ of degree $r + k_1 + 1$. By removing it, we obtain a graph $G'$, that is $(K_{1,r}; k_1)$-vertex stable. Obviously, $G$ is extremal if and only if $G'$ is extremal, which in turn takes place exactly when either $G' = G_1$ or $G' = G_2$. Thus

$$|E(G)| = \deg_G(v) + |E(G')| \geq (r + k_1 + 1) + |E(G(r, k_1))|$$
$$= (r + k_1 + 1) + \frac{1}{2}(k_1 + 1)(2r + k_1) = \frac{1}{2}(k_1 + 2)(2r + k_1 + 1)$$

with the lower bound being met only by $\hat{G}_1 \cong G(r, k_1 + 1)$ and $\hat{G}_2$ – the latter of which obtained from $K^{r+k_1-1}_{r+k_1+1}$ through addition of a total vertex.

**Case (3)** Let now $G$ be a $(K_{1,r}; k)$-vertex stable graph satisfying (3). Let $r$ be even and let
$$k_0 = \min\{k' : k' > (r-1)^2 - 2 \land 2 \nmid k'\}.$$

Consider first a $(K_{1,r}; k_0)$-vertex stable graph $G$. As the even value of $r$ implies that $(r-1)^2$ is odd, we may simplify the definition of $k_0$:
$$k_0 = \min\{k' : k' \geq (r-1)^2 \land 2 \nmid k'\} = (r-1)^2.$$

Suppose there exists a vertex $v$ in $G$ of degree $r + k_0$. By deleting it with incident edges, we obtain a $(K_{1,r}; k_0 - 1)$-vertex stable graph $G' = G - v$. Since $k_0 - 1 = k_1 + 1$, case (2) implies the following bound:

$$|E(G)| = \deg_G(v) + |E(G')| \geq r + k_0 + \frac{1}{2}(k_0 - 1 + 1)(2r + k_0 - 1) = \frac{1}{2}(k_0 + 1)(2r + k_0).$$

Suppose next $\Delta(G) \leq r + k_0 - 1$; by Lemma 2, $G$ is then $(r + k_0 - 1)$-regular and
$$|E(G)| = \frac{1}{2}(r + k_0 - 1)(r + k_0 + 1) = \frac{1}{2}((r + k_0)^2 - 1).$$

We may directly check, that
$$\frac{1}{2}((r+k_0)^2 - 1) < \frac{1}{2}(k_0+1)(2r+k_0),$$
therefore the only extremal graph is $K^{r+k_0-1}_{r+k_0+1}$.

It remains to show the inductive step. Assume the inductive hypothesis holds true for all odd $k'$ satisfying $k_0 < k' < k$; consider a $(K_{1,r};k)$-vertex stable graph $G$ and suppose there exists a vertex $v$ of degree $r+k$ in $G$. We may then remove $v$, thus getting a $(K_{1,r};k-1)$-vertex stable graph $G' = G - v$. Since $k-1$ is odd, there has to be a vertex $u$ of order $r+k-1$ in $G'$. By inductive hypothesis $G'' = G' - u = G - \{u,v\}$ is still $(K_{1,r};k-2)$-vertex stable and $|E(G'')| \geq \frac{1}{2}(r+k-3)(r+k-1)$. Therefore
$$|E(G)| = |E(G'')| + \deg_{G'}(v) + \deg_{G''}(u) \geq \frac{1}{2}(r+k-3)(r+k-1) + (r+k-1) + (r+k).$$

If, however, $\Delta(G) \leq r+k-1$, then $G$ is $(r+k-1)$-regular, whence
$$|E(G)| = \frac{1}{2}(r+k-1)(r+k+1).$$

We may directly check, that
$$\frac{1}{2}(r+k-1)(r+k+1) < \frac{1}{2}(r+k-3)(r+k-1) + (r+k-1) + (r+k),$$
equivalently
$$\frac{1}{2}((r+k)^2 - 1) < \frac{1}{2}((r+k)^2 - 1) + (r+k),$$
which concludes the proof of case (3).

**Case (4)** Let $r$ be even, $k_0$ be defined the same way as in case (3) and let $G$ be a $(K_{1,r};k)$-vertex stable graph. Consider now $(K_{1,r};k)$-vertex stable graphs for even values $k > k_0$.

We have already seen, that $\Delta(G) = r+k$. Let $v \in V(G)$ be such, that $\deg_G(v) = r+k$ and remove it from $G$ with incident edges. The resulting graph $G' = G - v$ satisfies the conditions of case (3); it is, therefore, $(r+k-2)$-regular. Hence
$$|E(G)| \geq \frac{1}{2}(r+k-2)(r+k) + (r+k) = \frac{1}{2}(r+k)^2.$$

We only have to show, that $G$ is in fact $(K_{1,r};k)$-vertex stable. Let $F$ denote the set of faults and $S$ the set of remaining vertices. Suppose $G$ is not $(K_{1,r};k)$-vertex stable. If $v \in S$, then $\deg_{G[S]}(v) = r$, since $\deg_G(v) = r+k$. Thus $v \in F$ and $|F \setminus \{v\}| = k-1$. However, $G'$ is $(K_{1,r};k-1)$-vertex stable, whence $\Delta(G[S]) \geq r$. Note, that we have also justified this way the fault tolerance of a graph considered in case (2), where the justification was omitted. The proof of the theorem is complete. □

## 4. Acknowledgements

The author would like to express his gratitude towards his Master's thesis supervisor Andrzej Żak for introducing the author to the problem and providing extensive support while dealing with it.

AGH University of Science and Technology, Faculty of Applied Mathematics, Al. A. Mickiewicza 30, 30-059 Kraków, Poland

*Email address*: `kuznar@agh.edu.pl, corresponding author`